\documentclass[12pt]{article}
\usepackage{amsmath}
\usepackage{amsfonts}
\usepackage{amssymb}
\usepackage{amsthm,amsmath}
\usepackage{amssymb,latexsym}
\usepackage{amscd,hyperref}
\usepackage[noend]{algpseudocode}
\usepackage{algorithm}
\usepackage{diagbox}
\usepackage{enumitem}
\usepackage{color}
\usepackage{graphicx}

\title{Stability Constrained Optimal Vessel Loading\\ \begin{small}
{Version 1.0}
\end{small}}
\author{Evangelos F. Magirou\footnote{Professor Emeritus, Athens University of Economics and Business. Email: efm@aueb.gr} }
\date{\today}

\begin{document}
\maketitle
		
\begin{abstract}
	We examine the revenue maximizing loading for the single voyage of a cargo vessel where metacentric stability is taken into account.  We formulate the problem as a maximization of a linear function with two linear  and one quadratic constraint.  The quadratic form in the constraint can be definite if heavier loads are to be placed low,  but not otherwise.  We consider a vessel similar to a 3500 TEU container carrier and perform several calculations using a Generalized Reduced Gradient commercial solver. 
	\end{abstract}

\section{Introduction}
We consider a vessel whose operator can select at will quantities from several types of cargoes in order to maximize revenue for a single voyage.  Each  type is characterized  by its density and its freight rate.  Loading constraints that have been taken into account in standard treatments \cite{Evans} are volume and weight ones, and specific rules are given for cargo selection.  In case volumes are additive we obviously have a linear programming problem \cite{MagPsar}.  For a container vessel we require integral solutions and a multidimensional integer knapsack problem results.

Stability questions can be easily incorporated in the above framework. We will show next that a simple stability constraint introduces an additional inequality constraint which is quadratic in the quantities loaded.  This type of problem (linear with one quadratic constraint) can be solved by a finite step algorithm in case the nonlinear constraint is convex i.e. with a positive semidefinite matrix. In our case, the quadratic form will indeed be positive definite if denser cargoes are to be placed low, but not otherwise.  

In the rest of this paper we formulate the stability constrained optimal loading problem, examine its properties and present some numerical calculations for a representative vessel.

\section{Problem Formulation}
 We assume that the available cargo types are indexed by $i$, have densities $d_i$, and freight rates $p_i$ monetary units per ton.  The vessel has a loading capacity of $C$ tons (Deadweight), while the volume of  loaded cargo can not exceed a maximum quantity $V$.  Assuming that the operator can select any nonnegative rational amount $x_i$ tons from the corresponding  type and wants to maximize total revenue, the problem is a linear programming one \cite{MagPsar}, provided we assume additivity in loaded weight and volume:
 
\begin{equation}
 \begin{array}{lll}
&\max_{x_i}\sum_{i=1}^{n}p_ix_i&\\
\text{such that }&&\\
&\sum_{i=1}^{n}x_i\le C&\\
&\sum_{i=1}^{n}x_i/d_i\le V&\\
&x_i\ge0&i=1,..,n\\
\end{array}
\end{equation}

\begin{figure}[ht]
	\centering
	\includegraphics[width=8cm, height=10cm]{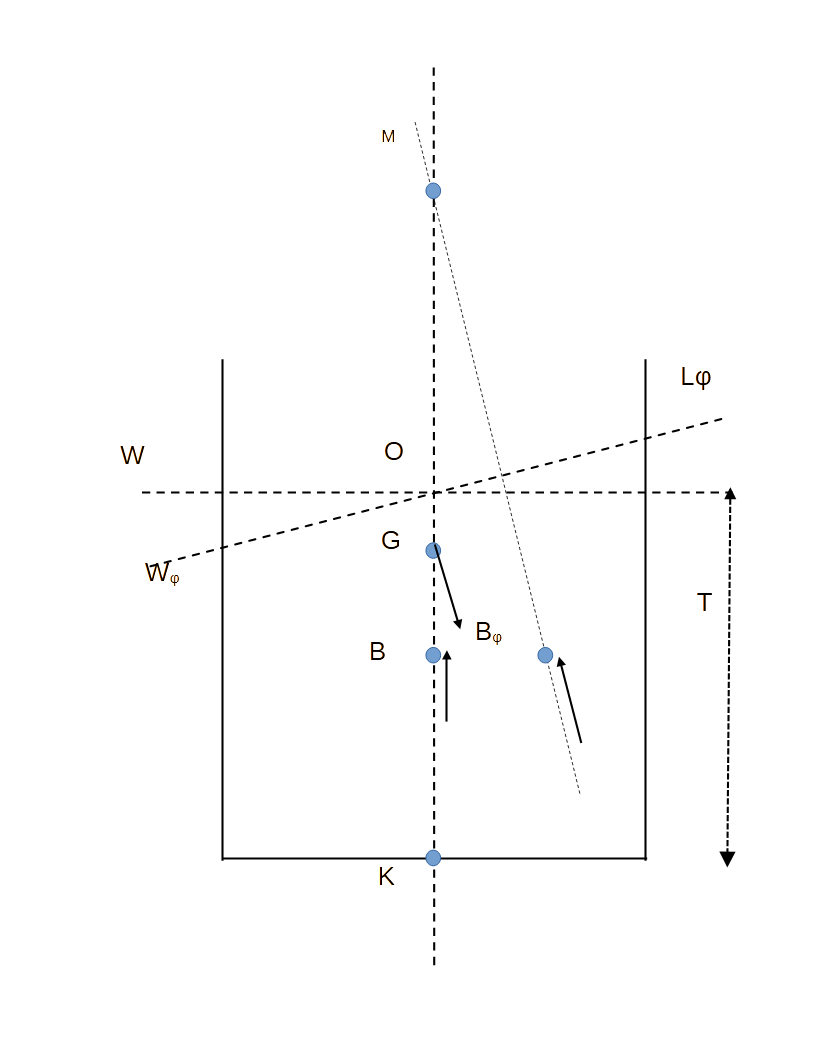}
	\caption{Metacenter Definition }
	\label{Meta}
	\end{figure}

The simplest stability consideration is expressed in terms of the (lateral) metacenter.  This is a particular feature of the vessel's shape and is the point through which all buoyancy force vectors pass (for small lateral inclinations). See  Figure \ref{Meta} and the textbook by Biran \cite{Biran}.   

For the ship to be stable the force couple  formed by  weight and  buoyancy must be equilibrating and this means that the  metacenter  should be above the vessel's center of weight by a comfortable margin.  The distance between  metacenter and  center of mass is the \textit{metacentric height}\footnote{In  Figure \ref{Meta} : \textit{K:Keel, M: Metacenter, WL:Waterline W$_\phi$L$_\phi$: Waterline for inclination $\phi$, B:Center of Buoyancy, G:Center of Mass, GM: Metacentric height.}} - $GM$ in Figure \ref{Meta}  and is proportional to the stabilizing moment; therefore for adequate stability the metacentric height should be sufficiently large. 

\begin{figure}[ht]	
	\centering
	\includegraphics[width=7cm, height=8cm]{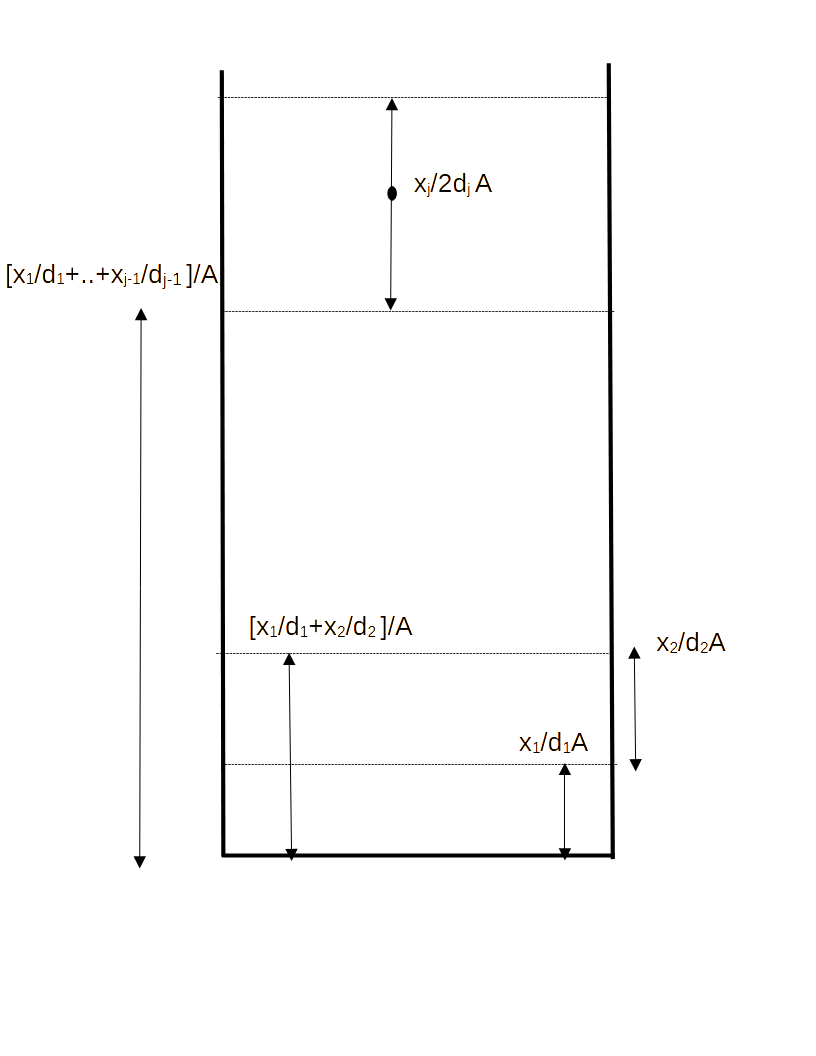}
	\caption{Vessel Loading\quad \textit{A:Area, $x_i,d_i$:Load and Density of freight $i$}}
	\label{Stack}
	\end{figure}

The metacentric height constraint places a restriction on the vessel's  center of mass whose location for a box shaped vessel is a quadratic function in the quantities loaded. This can be derived as follows:  Without loss of generality assume that the $j$-th cargo type will be loaded at the $j$-th position from the bottom.  Then cargo type  $j$ is placed  on top of $j-1$ others, which is a total height from the keel $(x_1/d_1+x_2/d_2+...+x_{j-1}/d_{j-1})/A$,  $A$ being the vessel's longitudinal area. The center of weight of cargo $j$ is a distance $x_j/(2Ad_j)$ above the lower ones, see Figure \ref{Stack}. 

If we let $M_{Vesl}$ be the weight of the empty vessel and $c_V$ the distance to keel of the empty vessel's center of mass, the center of mass  of the loaded vessel $y_{Mass}=KG$ in Figure \ref{Meta} is given by the relation
\begin{equation}\label{COM1}
	y_{Mass}[\sum_{i=1}^{n}x_i+M_{Vesl}]=c_V M_{vesl}+\frac{1}{A}\sum_{j=1}^{n}(x_1/d_1+x_2/d_2+..+x_{j-1}/d_{j-1}+x_j/2d_j)x_j
\end{equation}  
 
We can express the quadratic term using the matrix $W=\{w_{ij}\}=\{1/d_{\min(i,j)}\}$, i.e.
$$W=	\begin{bmatrix}
	d^{-1}_1       & d^{-1}_1 & d^{-1}_1& \dots & d^{-1}_1 \\
d^{-1}_1       & d^{-1}_2 &  d^{-1}_2 & \dots &  d^{-1}_2 \\
d^{-1}_1       & d^{-1}_2 &  d^{-1}_3 & \dots &  d^{-1}_3 \\
	\hdotsfor{5} \\
d^{-1}_1     &  d^{-1}_2 &  d^{-1}_3& \dots &  d^{-1}_n
	\end{bmatrix} $$

If $x'=(x_1,..,x_n)$ is the row vector of the selected loading  then \begin{equation}
\frac{1}{A}\sum_{j=1}^{n}(x_1/d_1+x_2/d_2+..+x_{j-1}/d_{j-1}+x_j/(2d_j))x_j=\frac{1}{2A}x'Wx
\end{equation}

and  thus if $1$ is a column  of ones, \eqref{COM1} can be written as   
\begin{equation}\label{COM}
y_{Mass}=KG=\frac{c_V M_{ves}+\frac{1}{2A}x'Wx}{x'1+M_{ves}}
\end{equation}

The location of the lateral metacenter is a distance $BM=\frac{I_{xy}}{\tilde{V}}$ from the center of buoyancy (see \cite{Biran} and Figure \ref{Meta}), where $I_{xy}$ is the moment of inertia of the water plane and $\tilde{V}$ is the submerged volume.  For a box shaped vessel whose hull is a parallelepiped of length $L$, breadth (beam) $B$ and draft $T$ the moment of inertia of the waterplane is $\frac{B^3L}{12}$, the submerged volume is $\tilde{V}=BLT$ and thus 
$$BM=\frac{B^3L}{12BLT}=\frac{B^2}{12T}.$$  
The  center of buoyancy $B$ is at a distance  $KB=T/2$ from the keel and thus the metacenter's keel distance is  $$KM=\frac{B^2}{12T}+\frac{T}{2}.$$  The draft $T$ is related to the displacement by $\rho LBT=M_{ves}+ x'1$, $\rho$ being the water density and thus 

\begin{equation}\label{Draft}
T=\frac{M_{ves}+ x'1}{\rho LB.}
\end{equation}

Using this  expression for the draft, the metacenter location $KM$ becomes 
$$KM=\frac{B^2}{12T}+\frac{T}{2}=\frac{1}{12}\frac{B^3L\rho}{M_{Ves}+x'1}+\frac{M_{Ves}+x'1}{2\rho LB}.$$

As stated earlier, a stability requirement is that the metacentric height is greater than a specified quantity say $\mu$,  i.e. $$y_{Mass}+\mu\le KM=\frac{B^2}{12T}+\frac{T}{2}.$$  Hence for our simple vessel the stability constraint becomes 
$$\frac{c_V M_{Ves}+\frac{1}{2A}xWx'}{x'1+M_{Ves}}+\mu \le \frac{B^2}{12T}+\frac{T}{2}.$$ 

Multiplying by $x'1+M_{Ves}$ and rearranging we obtain 

\begin{equation}\label{Stabil1}
\frac{1}{2BL}x'[W-\frac{1}{\rho}11']x+(\mu-\frac{M_{Ves}}{\rho BL})x'1\le \frac{{M_{Ves}^2}}{2\rho BL}+\frac{\rho B^3L}{12}-(\mu+c_V)M_{Ves}
\end{equation}

Thus the problem of maximizing the voyage revenue with a feasible loading is linear with a single quadratic constraint.

\begin{equation} \label{Main}
\begin{array}{lll}
&\max_{x_i}\sum_{i=1}^{n}p_ix_i&\\
\text{such that }&&\\
&\sum_{i=1}^{n}x_i\le C&\\
&\sum_{i=1}^{n}x_i/d_i\le V&\\
&\frac{1}{2BL}x'[W-\frac{1}{\rho}11']x+(\mu-\frac{M_{Ves}}{\rho BL})x'1\le \frac{{M_{Ves}^2}}{2\rho BL}+\frac{\rho B^3L}{12}-(\mu+c_V)M_{Ves}&\\
&x_i\ge0\qquad \qquad i=1,..,n&\\
\end{array}
\end{equation}

This is an linear optimization problem with a single quadratic constraint.  Several algorithms \cite{Panne}, \cite{Martein} exist which terminate in a finite number of iterations if  the constraints are convex.  This requires the  matrix $W-\frac{1}{\rho}11'$ to be positive semidefinite, but this does not always hold in our application. To examine this we use the following lemma\\

\textbf{Lemma}.
	\textit{Consider the $n$ real numbers  $(m_1,m_2,..,m_n)$ and   the symmetric $n\times n$ matrix $Q=\{q_{ij}\}$ with $q_{ij}=m_{min(i,j)}$  namely 
	$$Q=\begin{bmatrix}
		m_1&m_1&m_1&\dots &a_1\\
		m_1&m_2&m_2& \dots &  m_2 \\
		m_1&m_2&m_3& \dots &m_3\\
		\hdotsfor{5} \\
		m_1&m_2&m_3& \dots &m_n
	\end{bmatrix}.$$  Then $Q$ is congruent to a diagonal matrix with elements  $$(m_1,\;m_2-m_1,\ldots,\;m_k-m_{k-1},\ldots,\;m_n-m_{n-1}).$$ }
\begin{proof}
Using the standard method to obtain a congruent diagonal matrix (\cite{Fried} Sec. 6.8), we have the following reduction:  If 
$$
B=	\begin{bmatrix}
	1& 0&0& \dots &0\\
	1&1&0& \dots &0 \\
	1& 1&1& \dots &0\\
	\hdotsfor{5} \\
	1&1&1&\dots &1
\end{bmatrix}
$$
and 
$$D_Q=	\begin{bmatrix}
	m_1& 0&0& \dots &0\\
	0&m_2-m_1&0& \dots &0 \\
	0&0&m_3-m_2&\dots &0\\
	\hdotsfor{5} \\
	0&0&0&\dots &m_n-m_{n-1}
\end{bmatrix}
$$
We can verify by direct calculation that $$Q=BD_QB^T$$ 
\end{proof}
 
we can summarize the loading results in the following:\\

\textbf{Proposition}.
Consider the quadratic constraint matrix $W-\frac{1}{\rho}11'$ in \eqref{Main}. Then:
\textit{\begin{enumerate}
		\item The constraint matrix can not be positive semidefinite if some cargo is placed above a lower density one. 
		\item The constraint matrix is positive semidefinite provided dense cargoes are placed under lighter ones and furthermore the first cargo is not denser than water. 	
		\item The constraint matrix is negative semidefinite provided all cargoes are placed above lower density ones and  the first cargo is heavier than water.
	\end{enumerate}}
\begin{proof} 
Applying the Lemma to the matrix  $[W-\frac{1}{\rho}11']=\{a_{ij}\}=\{w_{ij}-\rho^{-1}\}$ we observe that it is congruent to a diagonal matrix with elements $d_1^{-1}-\rho^{-1},\;d^{-1}_2-d^{-1}_1,\;d^{-1}_3-d^{-1}_2,\ldots,\;d^{-1}_n-d^{-1}_{n-1}$.  Thus in Case 1 there exists at least one negative element in the diagonal matrix.  In Case 2 all diagonal elements are positive since density is decreasing while the first is also nonnegative since $d_1\le\rho$. Note that if  $d_1> \rho$ the matrix is not definite although it has only one negative element on the diagonal.  In Case 3 it is clear that all diagonal elements are negative.
\end{proof}	

In interpreting   Cases 2 and 3 above we can think of water acting as a ballast cargo, so if all cargoes and ballast  are loaded in decreasing density the constraint matrix  is positive semidefinite and problem \eqref{Main} has convex constraints.  If the cargoes and ballast are loaded in increasing density the matrix is negative semidefinite and we have a concave quadratic constraint in addition to convex linear constraints. 

We note that for cargoes denser than water the corresponding elements of the $W-1'1/\rho$ are negative and hence increasing them would increase metacentric height.  On the other hand the constraint matrix ceases to be positive. 	

In view of the Proposition it is  of practical interest to develop efficient algorithms for a non definite quadratic constraint.  We will not examine this further, and restrict our calculations to those that can be solved by the Generalized Reduced Gradient (GRG) nonlinear solver \cite{Lasdon} implemented in Excel.  

\section{An example} 

We present some calculations for a vessel similar to a 3500 TEU Container Carrier as in \textit{The Containership Register 2011} by Clarkson Research\footnote{Available online at https://www.aapa-ports.org/files/PDFs/CONTAINER\%20SHIP\%20SAMPLING.pdf}.  The vessel parameters are given in Table \ref{Tab:1}. For such a vessel, deadweight is about 45000 tons,  it can load about 3000 20 ft containers, each of volume 40 $m^3$ and thus 120000 $m^3$ volume capacity. Length is about 200 m, Beam 25-30 m and Weight of the empty vessel 15000 tons.  The center of mass position of the empty vessel is taken arbitrarily as 2 meters. 

\begin{table}[h]
	\centering
\begin{tabular}{|l|c|c|}
	\hline 
	Parameter &Units&Value  \\ 
	\hline
	Volume capacity &m$^3$&120000\\ \hline
	Deadweight&t&45000\\ \hline
	Beam&m&25\\ \hline
	Length&m&200\\ \hline
	Vessel Mass&t&15000\\ \hline
	Center of Mass of empty Vessel&m&2\\ 
	\hline 
\end{tabular}
\caption{Vessel Parameters \quad t: tons, m: meters}
\label{Tab:1}
\end{table}
 	
The cargo types appear in Table \ref{Tab:2} and are chosen so that they give interesting numerical results. For cargo types $i=1,2,3,4$  the  table gives the freight rates in monetary units per metric ton and density in tons per cubic meters.  Water density $\rho$ is assumed to be one, thus cargo densities are relative to water.  For convenience, indexing is in decreasing density $d_j>d_{j+1}$. We also assume there is a zeroth type not appearing in the table with unit density and zero freight rate, representing ballast. It is conceivable that some ballast might be loaded despite its zero freight rate in order to improve the vessel's handling, but we did not encounter such a case in our calculations. 

	\begin{table}[h]
		\centering
		\begin{tabular}{|l||c|c|c|c|}
		\hline 
		Type&1&2&3&4\\ 
		\hline
		Density&0.80&0.60&0.50&0.45\\
		Freight&4.50&5.00&5.10&5.50\\  
		\hline 
		\end{tabular}
		\caption{Cargo types}
		\label{Tab:2}
	\end{table}

We consider two cases for the loading order.  \textbf{Normal Loading}  occurs when all cargo types are placed under lower density ones and are lighter than water; thus the constraint matrix is positive definite by the Proposition. In \textbf{Reverse Loading} heavier cargoes are placed over lighter ones. As to stability requirements, we examine a \textbf{Relaxed Stability Margin} with metacentric height at least $\mu=4$ and a \textbf{Strict Stability Margin} with $\mu=6$.  We then solved the resulting four cases using the GRG  Excel Solver, and  then checked that the KT conditions are satisfied.  In the cases where  the algorithm  converged to a solution not satisfying the necessary conditions we restarted it from a different initial point until we obtained a solution satisfying them.  For Reverse Loading (Cases 1a, 2a) we encountered solutions that would satisfy the KT conditions but were not global optima. As expected this did not occur for Normal Loading where the constraints are convex.  

	\begin{table}[h]
	\centering
	\begin{small}
	\begin{tabular}{|c|c|c|c|c|c|c|c|c|c|c|c|}
		\hline 
		Case&\multicolumn{4}{c}{Loading - order}&&Revenue&$\mu$&Metac/ic&Volume&Metac.&Deadw.\\
		no.&\multicolumn{4}{c}{in $10^3$ tons}&&$10^3$&m&height m&$10^3$t&Multipl&Multipl\\
		\hline
		&1&2&3&4&Total&\multicolumn{5}{c}{}&\\ \hline
		1&8.5&9.1&&27.4&45.0&234.5&4&10.340&86.7&0.164&3.968\\ \hline
		2&33.2&1.7&&5.0&39.9&185.5&4&10.234&55.5&0.901&0.000\\ \hline
		&4&3&2&1&Total&\multicolumn{5}{c}{}&\\ \hline
		1a&2.7&&42.3&&45.0&226.3&6&10.340&76.5&0.106&4.225\\ \hline
		2a&&&&40.6&40.6&182.6&6&10.243&50.7&0.895&0.000 \\ \hline
	\end{tabular}
	\caption{Numerical Results}
	\label{Tab:3}
\end{small}
\end{table}

We present the  results of our calculations in Table \ref{Tab:3}. In all cases the quadratic constraint is binding.  If it were not, the solution would be a linear programming one with only two nonzero cargo types, as many as the Volume and Deadweight constraints.  This does not hold for Normal Loading with binding metacentric constraint where almost all types are loaded, as in Cases 1 and 2.  Fewer types are loaded in the Reverse Loading cases 1a and 2a, which seems reasonable since solutions will likely be on the edges of the linear constraint polyhedron; we did not examine the issue further. The revenue differences between Normal and Reverse loading are small, probably due to the small differences in the freight rates.  

Drastically increasing the stability requirement $\mu$ by 50\% has a dramatic effect: First, revenue decreases by about 20\% and total loaded cargo by more than 10\%. In the case of Reverse Loading the most profitable cargo type 4 is not to be used because it has to be placed at the bottom in spite of being the lightest.  Conversely, the least profitable type 1 is exclusively used in Case 2a since any other  cargo loaded has to be placed under it, thus seriously decreasing the metacentric height. Excluding type 1 is not advantageous since loading even limited quantities of the other types violates the metacentric constraint.   

We also present the constraint multipliers for all cases .  The volume multiplier is null since the constraint is not effective. The metacentric constraint multiplier is interesting in its own right but is also useful in order to assess the impact of the safety margin on revenue.  This can not  be directly obtained from multiplier information; however  an  examination of the quadratic constraint shows that assuming that the impact of the stability margin on the loading $x$ is small, the multiplier of $\mu$ is the metacentric constraint multiplier times the ship's displacement (cargo plus vessel's weight).  This can be verified by direct calculation.  In Case 1, if $\mu$ increases by 2.5\% from 4 to 4.1 the revenue  decreases by 1001 monetary units, or $1/235=0.4\%$. On the other hand, the sensitivity estimate is $\Delta \mu*0.164*(45000+15000)=984$ which differs from the calculated value of 1001 by only about $2\%.$  Note also that the deadweight multiplier is slightly smaller than the freight rates, reflecting the effect of the stability constraint.


\begin{thebibliography}{9}
\bibitem{Biran}   Biran, A.B. \textit{Ship Hydrostatics and 	Stability}. 1st Edition, Butterworth-Heinemann, An imprint of Elsevier, 2003.
\bibitem{Evans} Evans, J. and Marlowe, P. \textit{Quantitative Methods in Maritime Economics.}  1st Edition, Fairplay Publications, 1985.
\bibitem{Fried} Friedberg S., A. Insel  and L. Spence. \textit{Linear Algebra}. 4th Edition. Pearson Education, Inc., Upper Saddle River, New Jersey, 2003.
\bibitem{Lasdon} Lasdon L., R. Fox and M. Ratner: Nonlinear optimization using the generalized reduced gradient method, Revue Française d’Automatique, Informatique, Recherche Opérationnelle. Recherche Opérationnelle, tome 8, no V3, p. 73-103, 1974.
\bibitem{MagPsar}  	Magirou, E. and Psaraftis, H. and Christodoulakis, N. 	\textit{Quantitative Methods in Shipping: A Survey of Current Use and Future Trends}. Center for Economics Research Report E115 , Athens University of Economics and Business, 1992.  \textit{Available in the author's  Research Gate entry} 
\bibitem{Martein} Martein L. and S. Schaible:  On Solving a Linear Program with one Quadratic Constraint.  Rivisla di matematica per le scienze economiche e sociali - Anno 10", Fascicolo 1 *-2*, 1988.
\bibitem{Panne} van de Panne C: Programming with a Quadratic Constraint.  Management Science,  Vol. 12, Issue: 11, pp 798-815, 1966.
\end{thebibliography}
\end{document}